\documentclass[english]{article}
\usepackage[T1]{fontenc}
\usepackage[latin9]{inputenc}
\usepackage{calc}
\usepackage{amsmath}
\usepackage{amsthm}
\usepackage{amssymb}
\usepackage{stmaryrd}
\usepackage{graphicx}
\usepackage[numbers]{natbib}

\makeatletter

\usepackage{pdfpages}

\makeatother

\usepackage{babel}
\begin{document}
\noindent\fbox{\begin{minipage}[t]{1\columnwidth - 2\fboxsep - 2\fboxrule}%
In ``Contextuality in Random Variables: A Systematic Introduction''

by E. N. Dzhafarov, J. V. Kujala, and V. H. Cervantes

Cambridge, UK: Cambridge University Press, 2026

DOI 10.1017/9781009742221%
\end{minipage}}

\section*{Preface}

\subsection*{Subject of the book}

This book is about \emph{systems of random variables}.\index{system of random variables}
A system of random variables is a set of random variables ordered
in two ways: by their \emph{contents\index{content}} (the questions
the variables answer) and by their \emph{contexts\index{context}}
(conditions under which they are recorded). Such a system can be \emph{contextual}\index{system of random variables!contextual}
or \emph{noncontextual}.\index{system of random variables!noncontextual}
These terms, being derived from the word \emph{context}, have cognates
in many areas, but most of them are too syncretic to be helpful. They
do not distinguish contextuality from a broader notion of \emph{context-dependence}\index{context-dependence}.
The meaning of contextuality in this book is closest to that in quantum
physics, where its special forms are known as \emph{nonlocality} and
the \emph{Kochen--Specker contextuality}\index{Kochen--Specker!contextuality}.
Another related notion is that of \emph{selectiveness of influences}
in cognitive sciences, which is mathematically equivalent to noncontextuality
in the Kochen--Specker sense. However, our treatment is more general
and more abstract, applying to all systems of random variables, whether
they describe phenomena in quantum physics, computer science, or psychology.

The mathematical essence of contextuality is in the \emph{similarity\index{random variable!similarity}}
of random variables answering the same question (having the same content)
in different contexts. There is a principled way of measuring how
similar two such variables are: by computing the maximal probability
with which they could coincide \emph{if} they were jointly distributed
(de facto no two variables are if they are in different contexts).
A system is contextual if these same-content variables, in order to
be compatible with other variables in their respective contexts, have
to be more dissimilar than they are when considered in isolation,
outside their contexts. If this sounds less than clear, this book
will spell it out in great detail.

One fact, however, can be mentioned right away: The difference in
the similarity of two random variables considered within and without
their contexts is not due to any physical action of contexts upon
the variables. Considerations of causality may be involved in explaining
why the variables within this or that context have the joint distributions
they are observed or predicted to have. Contextuality, however, neither
requires nor admits explanations in terms of physical causes and effects.
It is a mathematical property, based on the fact that a random variable
is a function whose domain is determined by all other variables it
is jointly distributed with. Therefore a variable considered within
its context and a variable considered in isolation are simply two
different variables, even if they answer the same question. There
is no transition between the two, and they properly have to be denoted
by different symbols. In fact, contextuality can never be pinned to
specific variables, which it would have to if we dealt with causal
influences. When a system is contextual, we know that some of the
same-content variables are more dissimilar within their contexts than
they are without -- but one can never say what this difference in
similarity is for any particular pair. Whether a system is or is not
contextual is a \emph{system-wise} (and \emph{system-wide}) property
of the system. In particular, it is easy to construct examples of
contextual systems that become noncontextual (or vice versa) following
a modification of any single variable's distribution in the system.
This applies, in particular, to variables recorded prior to all other
variables in their contexts, excluding thereby any possibility of
causal influences being involved.

\subsection*{Philosophy of the book}

It is the leitmotif of this book that contextuality-related concepts
are purely mathematical constructs. To preclude a hasty reaction to
this statement: Something being purely mathematical does not prevent
it from being useful and widely applicable. Matrices, integrals, and
arithmetic operations are purely mathematical constructs, based on
no empirical assumptions and subject to no empirical verification.
They nevertheless have numerous applications, and they are highly
useful in all of them. It is simply that no physical considerations
are involved in determining if a matrix is invertible or if the value
of an integral is finite.

Although the interplay between the mathematical and the empirical
is a deep philosophical issue which is outside the reach of this book,
one can acknowledge some obvious facts. Mathematical analysis (in
our case, contextuality analysis) always applies to mathematical descriptions
of what is being studied. These descriptions (in our case, systems
of random variables) are always non-unique, but they are not arbitrary.
A situation they describe constrains them and makes them interrelated.
One can say that it is an empirical property of a situation to afford
a certain class of mathematical descriptions and only them. However,
mathematical properties of different descriptions in this class (in
our case, whether they are contextual or not) may differ. Interpreting
contextuality-related concepts in terms of hypothetical natural phenomena
may be innocuous if these phenomena are viewed as mere analogies or
metaphors. But on a deeper level of analysis such interpretations,
as a rule, obscure understanding. Thus, as thoroughly explained in
this book, a \emph{hidden variable model\index{hidden variable model}}
with \emph{free choice\index{hidden variable model!free choice}}
but \emph{context-dependent mapping\index{context-dependence}} of
the hidden variables into observable ones is equivalent to a hidden
variable model with \emph{context-independent mapping} but compromised
freedom of choice. This simple mathematical fact looks paradoxical
if not impossible if one interprets hidden variables as physical states
or forces, and freedom of choice as a property of a psychological
(or machine-generated) process of decision making.

Our approach largely follows the great philosophical traditions of
logical positivism\index{logical positivism} and pragmatism\index{pragmatism}.
Although they are often eschewed nowadays by scientists and philosophers,
their non-radical versions serve as highly useful guides. Logical
positivism prevents us from constructing theories based on vague intuitions
and connotations of words. Pragmatism complements this by making us
realize that competing descriptions are not really different unless
they have distinct observable consequences. While these philosophical
traditions by no means exhaust all philosophy of science, they form
a necessary part of it. They are in fact indispensable in any truly
careful conceptual analysis. One can see this whenever one tries to
achieve sufficient clarity in dealing with such philosophy-laden issues
as \emph{fundamental stochasticity}\index{stochasticity!fundamental},
\emph{counterfactual definiteness}\index{counterfactual!definiteness},
or the nature of random variables. Contextuality turns out to be related
to several such issues. Can fundamental stochasticity be conceptually
distinguished from \emph{ensemble stochasticity}\index{stochasticity!ensemble},
and is there a way to demonstrate that one and not the other of these
forms of stochasticity is present in a given situation? It turns out
that contextuality in combination with another property, called \emph{non-disturbance},
can be a key to answering both these questions affirmatively. Is it
possible to assign truth values to counterfactual statements about
values of random variables, and if so, can counterfactual definiteness
(the statement that the value of a variable would be the same if it
were recorded in another context) be distinguished from contextuality?
Again, a careful conceptual analysis allows one to answer yes to both
these questions (with the clarification that it is not the same variable
in different contexts one is dealing with but different variables
with the same content). Random variables generate their values (say,
$Yes$ and $No$) with certain probabilities, but to corroborate this
one should be able to observe one and the same random variable repeatedly.
How does one know this is the same variable rather than a new variable
with the same distribution in every new observation? A careful analysis
of this conceptual problem leads one to the notion of \emph{probabilistic
couplings\index{coupling}} prominently represented in this book.

\subsection*{Genre of the book I: A textbook}

This book is almost entirely based on the authors' own published work,
which, however, the book cites less than sporadically. In general,
the book contains very few literature references, and those it does
are mostly of the historical nature. The book is written with a deliberate
intent \emph{not} to discuss the complex and plentiful literature
related to contextuality. This is not a slight to the important achievements
of the researchers in the contextuality field. It is simply that the
book is more of a textbook than a scholarly survey, and it follows
the style of its genre.

We should qualify this statement: In what sense is this book a textbook?
One sense in which it is, is that it presents the material systematically,
rigorously, and with no previous knowledge thereof required. It is
also comprehensive, but only in that it comprehensively covers the
basics of contextuality on an abstract mathematical level. This excludes
many specialized and domain-specific topics, those that may be of
interest in, or even central to, a substantive area of research, such
as computer science, psychology, or quantum physics in which most
of the contextuality research has been conducted.\footnote{The reader interested in how contextuality applies to quantum physics
can find a wealth of information in the available surveys: Yeong-Cherng
Liang, Robert W. Spekkens, \& Howard M. Wiseman, \emph{Physics Reports}
\textbf{506}, 1 (2011), and Costantino Budroni, Adán Cabello, Otfried
Gühne, Matthias Kleinmann, \& Jan-Åke Larsson, \emph{Review of Modern
Physics} \textbf{94}, 045007 (2022).}

Also, this book is a textbook in the sense of being accessible. It
is mostly confined to systems with a finite number of random variables,
and to variables with a finite number of possible values, mostly,
just two. This allows the reader to focus on all conceptual issues
without getting into technical details requiring previous knowledge
of measure-theoretic constructs (although these are introduced too,
as optional reading). The book is written to be accessible to college
students or even to advanced high-schoolers with only modest knowledge
of mathematics: The basic set-theoretic notions and notation, elementary
aspects of probability theory, and, in places, basics of linear algebra.
Each chapter is followed by exercises with complete solutions. They
provide additional examples and clarifications, and sometimes offer
additional theoretical material. Technical issues and involved proofs
are put in special sections marked by asterisks. They may be skipped,
although hopefully they will not be by better prepared and more skeptical
readers. Proofs are often conducted by means of easily generalizable
examples. The latter does not mean, however, that the book compromises
on deductive reasoning. Quite the opposite, deductive reasoning permeates
the entire book, so that the reader can always see how the issues
and results logically follow from a few initial definitions, principles,
and conventions.

\subsection*{Genre of the book II: A treatise for experts}

We should, however, mention the sense in which this book is \emph{not}
a textbook as commonly understood: It does not present a widely agreed
upon set of well-established topics and results. I am not sure such
a set exists. At the contemporary stage of contextuality studies any
attempt to present a systematic theory of contextuality which is not
domain-specific is bound to be an original approach. The present attempt
is based on the theory dubbed \emph{Contextuality-by-Default},\index{Contextuality-by-Default theory}
developed by Janne V. Kujala and myself since 2014. Other colleagues,
prominently including Víctor H. Cervantes, joined us in this development
at its later stages. One feature of the theory is that random variables
are systematically identified not only by their \emph{contents} but
also by their \emph{contexts}. This double-identification allows one
to stay well within the confines of classical probability theory,
or more precisely, classical theory of random variables. In fact,
the double-identification is a necessary requirement of this theory.
A systematic use of the language of random variables has mathematical
advantages over the more common probabilistic treatments of contextuality
in the language of events. This is because random variables add to
probability spaces the flexibility of functions mapping them into
each other. Thus, the assumptions traditionally considered necessary
for the derivation of the Bell-type inequalities, such as \emph{outcome
determinism} and \emph{factorizability\index{hidden variable model!factorizability}},
are satisfied automatically once the underlying hidden variables are
treated as true random variables. Our treatment of random variables,
however, puts a greater emphasis than is usual on the situations when
random variables do not possess joint distributions. As a consequence,
there is a greater emphasis on the theory of \emph{probabilistic couplings}\index{coupling}.
The book also attaches greater than usual importance to the systems
of \emph{dichotomous random variables}, and to the possibility of
redefining any system of variables as a system of dichotomous ones.
We even consider, without committing to it, the possibility of \emph{radical
dichotomism}\index{dichotomization!radical}, according to which the
contextuality status of a system is determinable only after the system
is dichotomized.

The greatest deviation from the traditional accounts of contextuality,
however, is in that we do not constrain the systems of random variables
by the requirement that the distribution of a variable only depends
on its content and not on its context. In other words, if two variables
answer the same question, we allow their distributions to differ.
The assumption that they must be the same is justifiable for some
quantum-theoretic systems, but it fails to apply to many others (e.g.,
those describing successive measurements, or the multiple-slit experiments).
Moreover, this constraint almost never holds outside physics, for
example, in behavioral and social sciences. Our book shows a principled
way of isolating and quantifying contextuality in all such situations.
(The reader who is only interested in \emph{undisturbed systems},
those with context-independent distributions, still might benefit
from reading this book, both because such systems are prominently
discussed in the book as special cases, and because it is always useful
to see a concept in a more general mathematical setting.) The book
also contains some new, previously unpublished or incompletely published
developments, such as the hierarchical contextuality measure and a
general theory of dichotomizations in structured spaces. Overall,
in spite of its technical accessibility, the material presented in
this book is advanced enough to be of interest to scholars from a
broad spectrum of disciplines: From mathematics to philosophy to quantum
physics to computer science to behavioral and social sciences.

\subsection*{Terminology of the book}\index{terminology}

There is no universally or even widely accepted terminology in the
contextuality literature. Even most basic notions (such as \emph{signaling})
have multiple terminological variants, and some frequently used terms
(such as \emph{realism}) allow for multiple interpretations. The terminology
adopted in this book therefore follows the internal logic of the theory
being developed rather than literature sources. We begin by calling
the set of random variables sharing a context a \emph{bunch} of variables
(because they form a single random variable), and we call the set
of random variables sharing a content a \emph{connection} of variables
(because they relate to each other probabilistically unrelated bunches).
Virtually all other terms that are not traditional are derivations
from these two: \emph{well-bunched}, \emph{well-connected}, \emph{consistently
connected}. Of course, like in most systematic treatments of mathematical
subjects, in ours one can find degrees of variation in the use of
standard mathematical terms. Thus, the standard term \emph{maximal
coupling} acquires its generalized version in \emph{multi-maximal
coupling}. We also introduce the term \emph{agglutinativity} to describe
the basic property of random variables: A set of variables is jointly
distributed if so are elements of a chain of its overlapping subsets.
I believe that with a bit of patience and goodwill the reader will
find our terminology thoughtfully designed and, due to its generality,
convenient to use.

\subsection*{Acknowledgements }

In the development of the theory underlying this book the authors
benefited from communications and debates with many colleagues, of
whom I would like to especially mention Samson Abramsky, Harald Atmanspacher,
Guido Bacciagaluppi, Acacio de Barros, Jerome Busemeyer, Adán Cabello,
Matthew Jones, Philippe Grangier, Andrei Khrennikov, and Pawe\l{}
Kurzy\'{n}ski. In the development of the precursor theory, that of
selective influences in human cognition, I would like to thank Hans
Colonius, R. Duncan Luce, Richard Schweickert, Patrick Suppes, Robin
Thomas, and James T. Townsend. I am sure I have forgotten some names
I should have mentioned.

\medskip{}

\begin{flushright}
Ehtibar N. Dzhafarov 
\par\end{flushright}

\begin{flushright}
Prague, Czech Republic 
\par\end{flushright}

\begin{flushright}
\includegraphics[scale=0.15]{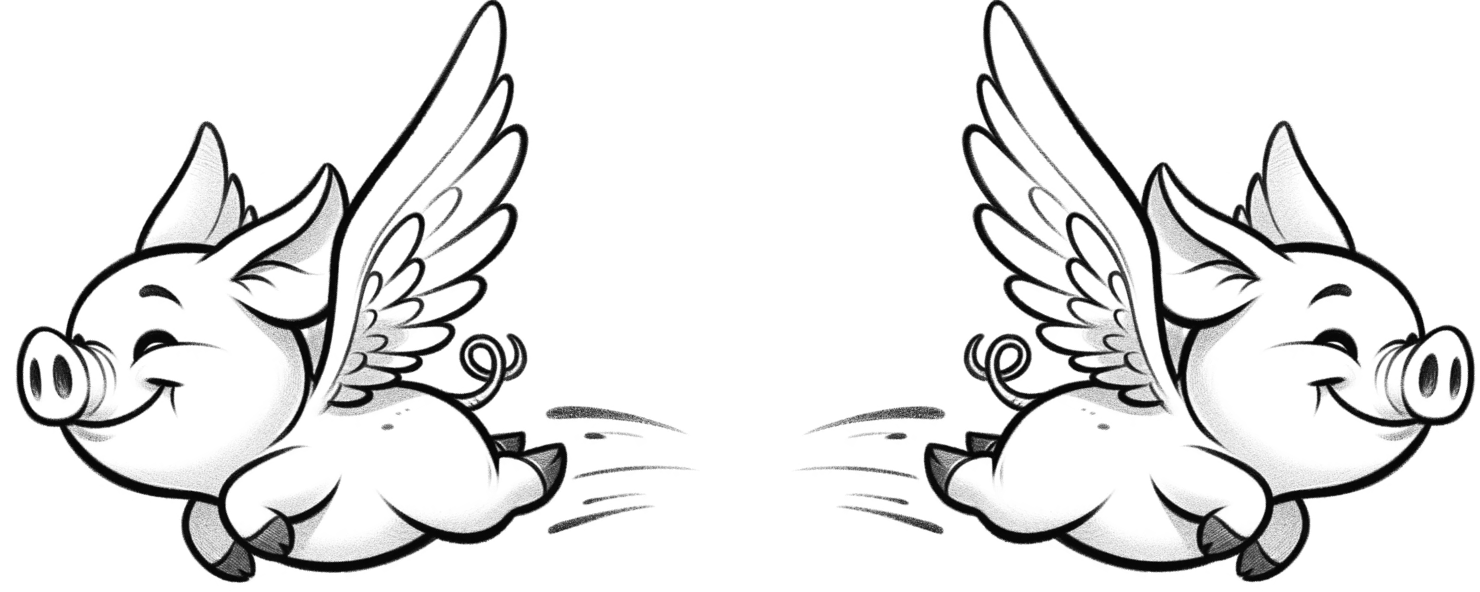} 
\par\end{flushright}

\newpage{}
\begin{flushleft}
\includegraphics[scale=0.7]{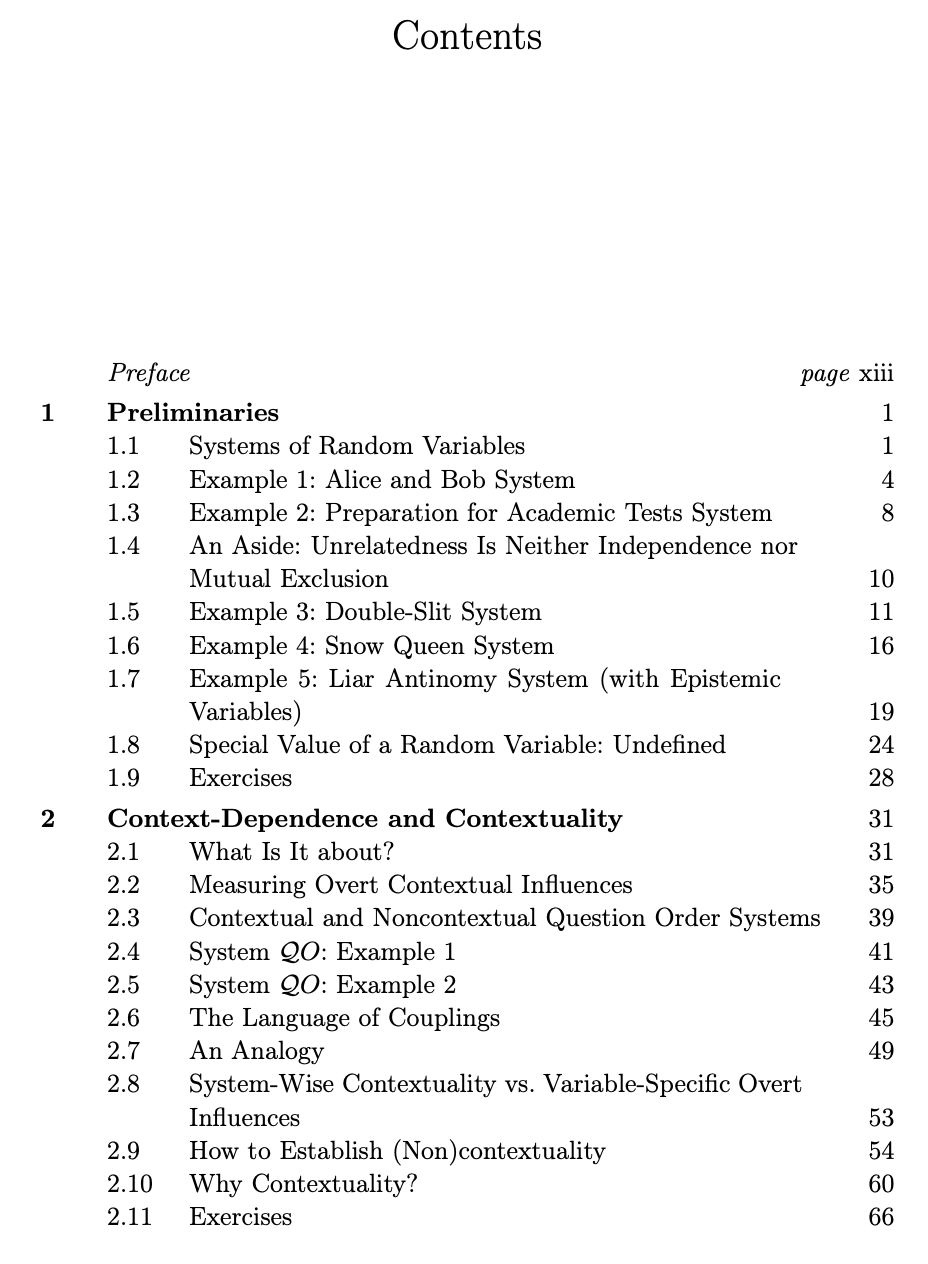}
\par\end{flushleft}

\begin{flushleft}
\includegraphics[scale=0.7]{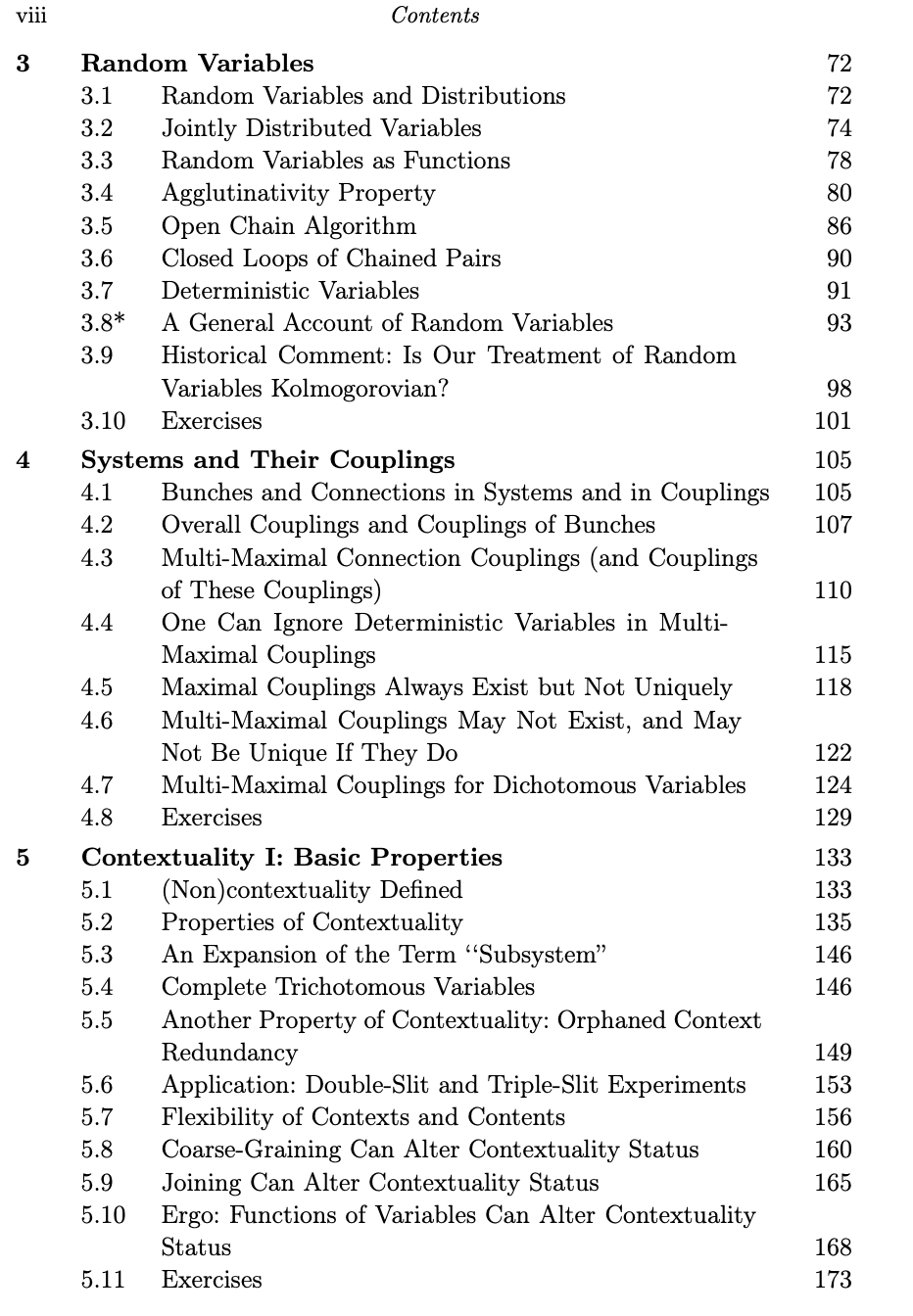}
\par\end{flushleft}

\includegraphics[scale=0.7]{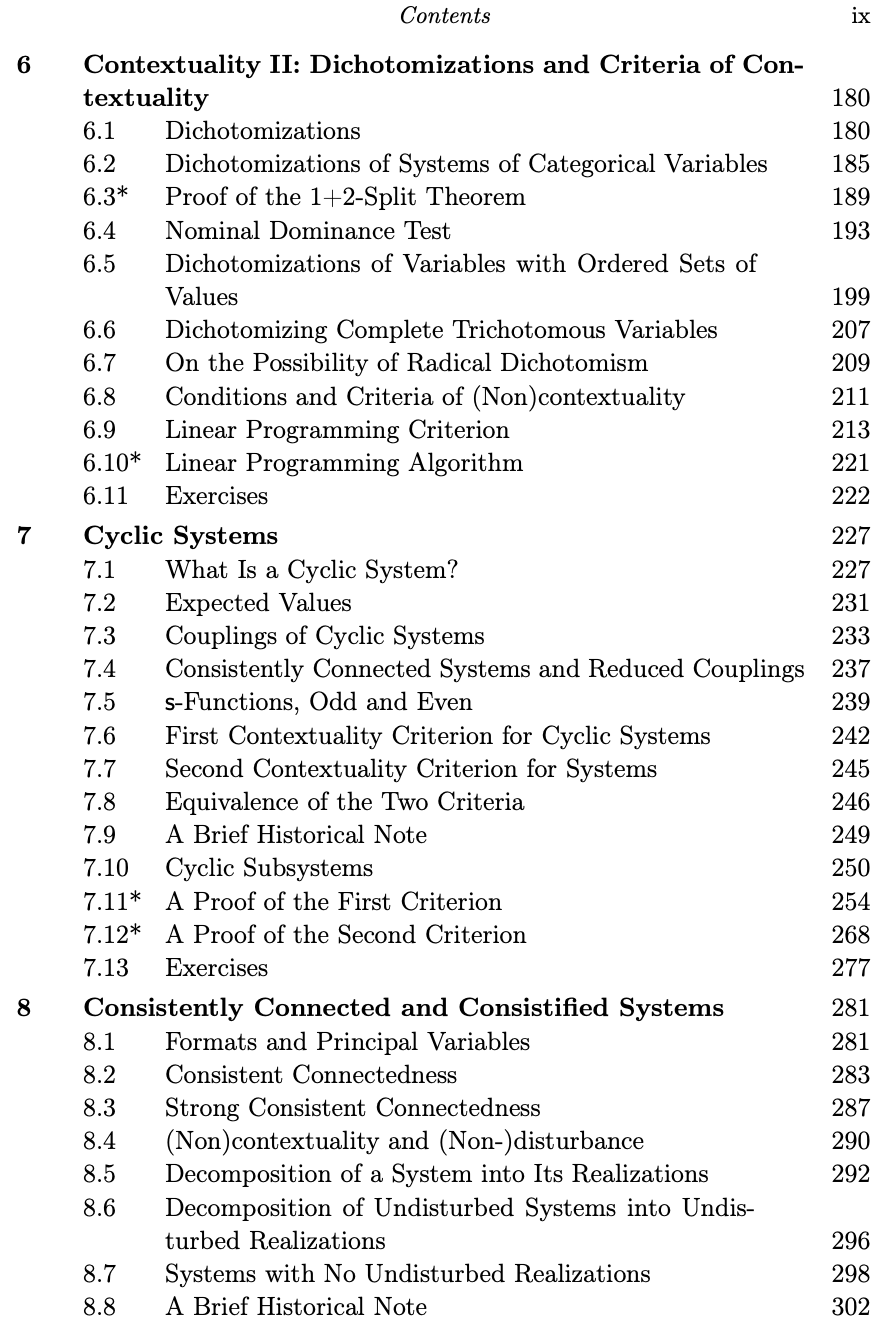}
\begin{flushleft}
\includegraphics[scale=0.7]{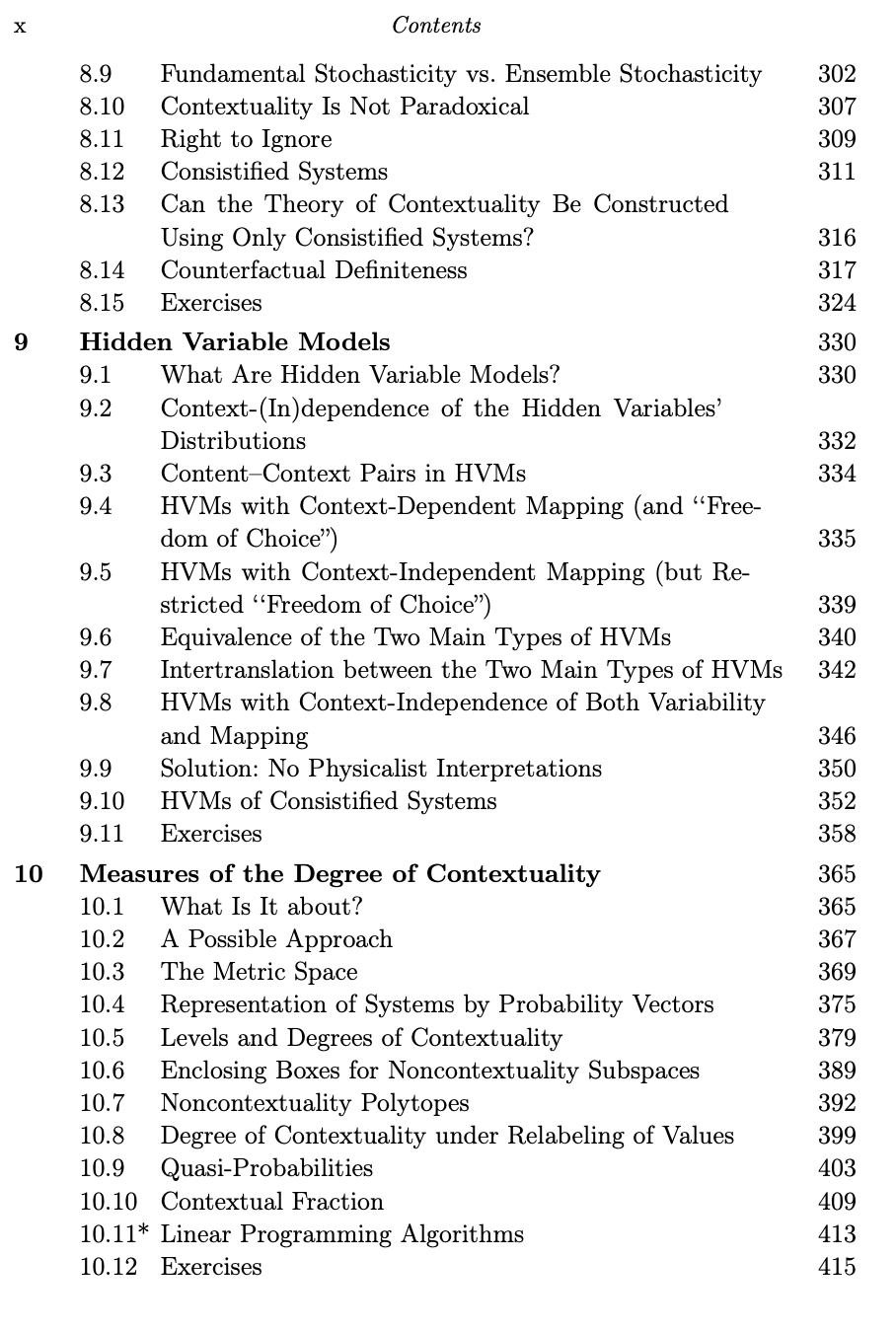}
\par\end{flushleft}

\begin{flushleft}
\includegraphics[scale=0.7]{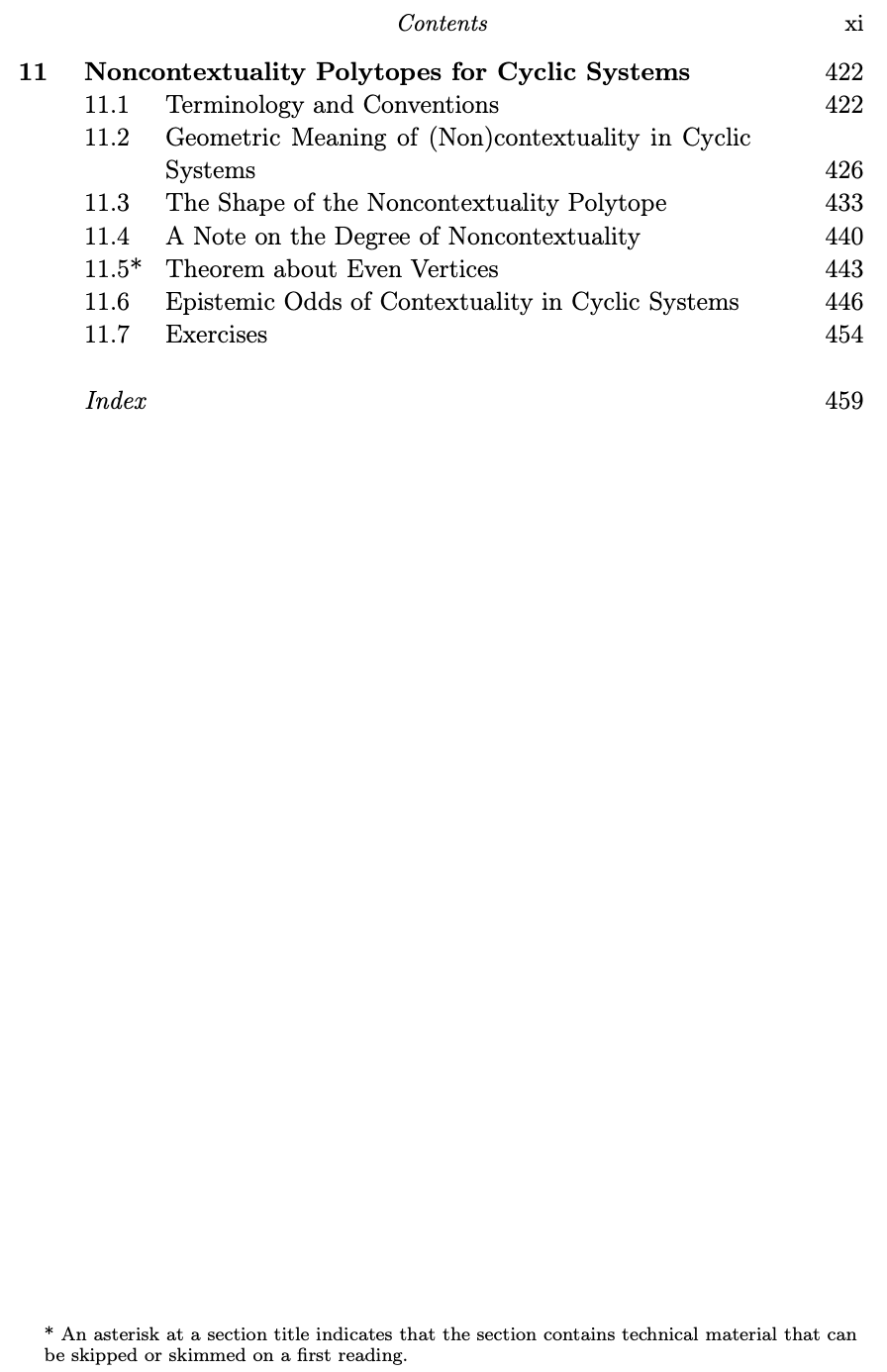}
\par\end{flushleft}
\end{document}